\documentclass[]{article}
\usepackage{amssymb}
\usepackage{amsmath}
\usepackage{latexsym,a4wide}

\input xypic
\input xy
\xyoption{all}

\newtheorem{example}{Example}[section]
\newtheorem{defn}[example]{Definition}
\newtheorem{prop}[example]{Proposition}
\newtheorem{thm}[example]{Theorem}

\newtheorem{ex}[example]{Example}
\newenvironment{pf}{\noindent \textbf{Proof:} }{ $\Box$ \mbox{} }
\newcommand{\llabto}[2]{\stackrel{#2}{\rule[0.52ex]{#1em}{0.099ex}\hspace{-0.4em}\longrightarrow}}
\def\leq{\leqslant}
\def\geq{\geqslant}
\date{}
\begin{document}

\author{Z. Arvasi and E. Ulualan}
\title{On Algebraic Models for Homotopy 3-Types}
\maketitle

\begin{abstract}
We explore the relations among quadratic modules, 2-crossed
modules, crossed squares and simplicial groups with Moore complex
of length $2$.
\end{abstract}

\section*{Introduction}
Crossed \ modules defined by Whitehead, \cite{wayted}, are
algebraic models of connected (weak homotopy)  2-types. Crossed
squares as introduced by Loday and Guin-Walery, \cite{walery},
model connected 3-types. Crossed $n$-cubes model connected
$(n+1)$-types, (cf. \cite{Porter}). Conduch\'{e}, \cite{Conduce},
gave an alternative model for connected 3-types in terms of
crossed modules of groups of length 2 which he calls `$2$-crossed
module'. Conduch\'{e} also constructed (in a letter to Brown in
1984) a $2$-crossed module from a crossed square. Baues,
\cite{baus2}, gave the notion of quadratic module which is a
$2$-crossed module with additional `nilpotency' conditions. A
quadratic module is thus a `nilpotent' algebraic model of
connected $3$-types. Another algebraic model of connected 3-types
is `braided regular crossed module' introduced by Brown and
Gilbert (cf. \cite{rb-ng}). These notions are then related to
simplicial groups. Conduch\'{e} has shown that the category of
simplicial groups with Moore complex of length $2$ is equivalent
to that of 2-crossed modules. Baues
gives a construction of a quadratic module from a simplicial group in \cite%
{baus2}. Berger, \cite{berger}, gave a link between $2$-crossed
modules and double loop spaces.

Some light on the $2$-crossed module structure was also shed by
Mutlu and Porter, \cite{mutpor}, who suggested ways of
generalising Conduch\'{e}'s construction to higher $n$-types. Also
Carrasco-Cegarra, \cite{cc}, gives a generalisation of the
Dold-Kan theorem to an equivalence between simplicial groups and a
non-Abelian chain complex with  a lot of extra structure,
generalising $2$-crossed modules.

The present article aims to show some relations among algebraic
models of connected $3$-types. Thus the main points of this paper
are:

(i) to give a complete description of the passage from a crossed
square to a 2-crossed module by using the `Artin-Mazur'
codiagonal functor and prove directly a 2-crossed module
structure;

(ii) to give a functor from 2-crossed modules to quadratic modules
based on Baues's work (cf. \cite{baus2});

(iii) to give a full description of a construction of a quadratic
module from a simplicial group by using the Peiffer pairing
operators;

(iv) to give a construction of a quadratic module from a crossed square.

Therefore, the results of  this paper can be summarized in the
following commutative diagram
\begin{equation*}
\xymatrix{ &\textbf{SimpGrp$_{\leq 2}$}\ar[dl]\ar[d]\ar[dr]&&\\
\textbf{QM}&\textbf{X$_2$Mod}\ar[l]\ar[u]& \textbf{Crs$^2$}
\ar@/^1.2pc/[ll] \ar[ul]\ar[l]&\\&&&}
\end{equation*}
where the diagram is commutative, linking the constructions given
below.

\noindent \textbf{Acknowledgements}. The authors wishes to thank
the referee for helpful comments and improvements to the paper.

\section{\label{s1} Preliminaries}

We refer the reader to May's book (cf. \cite{may}) and
Artin-Mazur's, \cite{Artin}, article for the basic properties of
simplicial groups, bisimplicial groups, etc.

A simplicial group $\mathbf{G}$ consists of a family of groups
${G_{n}}$ together with face and degeneracy maps
$d_{i}^{n}:G_{n}\rightarrow G_{n-1}$, $0\leq i\leq n$ $(n\neq 0)$
and $s_{i}^{n}:G_{n}\rightarrow G_{n+1}$, $0\leq i\leq n$
satisfying the usual simplicial identities given by May. In fact
it can be completely described as a functor $\mathbf{G}:\Delta
^{op}\rightarrow \mathbf{Grp}$ where $\Delta $ is the category of
finite ordinals.

Given a simplicial group $\mathbf{G}$, the Moore complex
$(\mathbf{NG}, \partial )$ of $\mathbf{G}$, is the (non-Abelian)
chain complex defined by;
\begin{equation*}
{NG}_{n}=\ker d_{0}^{n}\cap \ker d_{1}^{n}\cap \cdots \cap \ker d_{n-1}^{n}
\end{equation*}%
with $\partial _{n}:NG_{n}\rightarrow NG_{n-1}$ induced from $d_{n}^{n}$ by
restriction.

The $n^{th}$ \emph{homotopy group} $\pi _{n}(\mathbf{G})$ of $\mathbf{G}$ is
the $n^{th}$ homology of the Moore complex of $\mathbf{G}$, i.e.
\begin{equation*}
\pi _{n}(\mathbf{G})\cong H_{n}(\mathbf{NG},\partial
)=\left(\bigcap\limits_{i=0}^{n}\ker
d_{i}^{n}\right)/d_{n+1}^{n+1}\left( \bigcap\limits_{i=0}^{n}\ker
d_{i}^{n+1}\right) .
\end{equation*}
The Moore complex carries a lot of fine structure and this has
been studied, e.g. by Carrasco and Cegarra (cf. \cite{cc}), Mutlu
and Porter (cf. \cite{Mutlu, mutlu2, mutpor}).

Consider the product category $\Delta \times \Delta $ whose objects are pairs $%
([p],[q])$ and whose maps are pairs of weakly increasing maps. A
(contravariant) functor $\mathbf{G.,.}:(\Delta \times \Delta
)^{op}\rightarrow \mathbf{Grp}$ is called a bisimplicial group. Hence $%
\mathbf{G.,.}$ is equivalent to giving for each $(p,q)$ a group
$G_{p,q}$ and morphisms:
\begin{equation*}
\begin{tabular}{ll}
$d_{i}^{h}:G_{p,q}\rightarrow G_{p-1,q}$ &  \\
$s_{i}^{h}:G_{p,q}\rightarrow G_{p+1,q}$ & $0\leq i\leq p$ \\
$d_{j}^{v}:G_{p,q}\rightarrow G_{p,q-1}$ &  \\
$s_{j}^{v}:G_{p,q}\rightarrow G_{p,q+1}$ & $0\leq j\leq q$%
\end{tabular}%
\end{equation*}%
such that the maps $d_{i}^{h},s_{i}^{h}$ commute with $d_{j}^{v},s_{j}^{v}$
and that $d_{i}^{h},s_{i}^{h}$ (resp. $d_{j}^{v},s_{j}^{v}$) satisfy the
usual simplicial identities.

We think of $d_{j}^{v},s_{j}^{v}$ as the vertical operators and $%
d_{i}^{h},s_{i}^{h}$ as the horizontal operators. If $\mathbf{G.,.}$
is a bisimplicial group, it is convenient to think of an element of
$G_{p,q}$ as a product of a $p$-simplex and a $q$-simplex.

\section{\label{s2}2-Crossed Modules from Simplicial Groups}

Crossed modules were initially defined by Whitehead as models for
connected 2-types. As explained earlier, Conduch{\'{e}},
\cite{Conduce}, in 1984 described the notion of $2$-crossed module
as models for connected 3-types.

\emph{A crossed module} is a group homomorphism $\partial :M\rightarrow P$
together with an action of $P$ on $M$, written $^{p}m$ for $p\in P$ and $%
m\in M$, satisfying the conditions $\partial (^{p}m)=p\partial
(m)p^{-1}$ and $^{\partial m}m^{\prime }=mm^{\prime }m^{-1}$ for
all $m,m^{\prime }\in M,p\in P$. The last condition is called the
`Peiffer identity'.

The following definition of $2$-crossed module is equivalent to
that given by Conduch{\'{e}}.

\emph{A 2-crossed module} of groups consists of a complex of groups%
\begin{equation*}
\xymatrix@C=0.5cm{
  L \ar[rr]^{\partial_2} && M \ar[rr]^{\partial_1} && N  }
\end{equation*}%
together with (a) actions of $N$ on $M$ and $L$ so that $\partial
_{2},\partial _{1}$ are morphisms of $N$-groups,  and (b) an
$N$-equivariant function
\begin{equation*}
\{\quad ,\quad \}:M\times M\longrightarrow L
\end{equation*}%
called a Peiffer lifting. This data must  satisfy the following
axioms:
\begin{equation*}
\begin{array}{lrrll}
\mathbf{2CM1)} &  & \partial _{2}\{m,m^{\prime }\} & = & \left( ^{\partial
_{1}m}m^{\prime }\right) mm^{\prime }{}^{-1}m^{-1}\newline
\\
\mathbf{2CM2)} &  & \{\partial _{2}l,\partial _{2}l^{\prime }\} & = &
[l^{\prime },l]\newline
\\
\mathbf{2CM3)} &  & (i)\quad \{mm^{\prime },m^{\prime \prime }\} & = &
^{\partial _{1}m}\{m^{\prime },m^{\prime \prime }\}\{m,m^{\prime }m^{\prime
\prime }m^{\prime }{}^{-1}\}\newline
\\
&  & (ii)\quad \{m,m^{\prime }m^{\prime \prime }\} & = & \{m,m^{\prime
}\}^{mm^{\prime }m^{-1}}\{m,m^{\prime \prime }\}\newline
\\
\mathbf{2CM4)} &  & \{m,\partial _{2}l\}\{\partial _{2}l,m\} & = &
^{\partial _{1}m}ll^{-1}\newline
\\
\mathbf{2CM5)} &  & ^{n}\{m,m^{\prime }\} & = & \{^{n}m,^{n}m^{\prime }\}%
\newline
\end{array}%
\end{equation*}%
\newline
for all $l,l^{\prime }\in L$, $m,m^{\prime },m^{\prime \prime }\in M$ and $%
n\in N$.

Here we have used $^{m}l$ as a shorthand for $\{\partial_2l, m\}l$
in condition $\mathbf{2CM3)} (ii)$ where $l$ is $\{m,m^{\prime
\prime}\}$ and $m$ is $mm^{\prime}(m)^{-1}$. This gives a new
action of $M$ on $L$. Using this notation, we can split
$\mathbf{2CM4)}$ into two pieces, the first of which is
tautologous:
\begin{equation*}
\begin{array}{lrrll}
\mathbf{2CM4)} & \quad (a) \quad \{\partial_2 l,m\} & = &
^{m}l(l)^{-1},
\\
& \quad   (b) \quad \{m, \partial_2 l \}  & = & (^{\partial_1 m}
l)(^{m}l^{-1}).
\end{array}
\end{equation*}
The old action of $M$ on $L$, via $\partial_1$ and the $N$-action
on $L$, is in general distinct from this second action with
$\{m,\partial_2 l\}$ measuring the difference (by $\mathbf{2CM4)}
(b)$). An easy argument using $\mathbf{2CM2)}$ and $\mathbf{2CM4)}
(b)$ shows that with this action, $ ^{m} l$, of $M$ on $L$,
$(L,M,\partial_2)$ becomes a crossed module.

A morphism of $2$-crossed modules can be defined in an obvious
way. We thus define the category of $2$-crossed modules denoting
it by \textbf{X}$_{2}$\textbf{Mod}.

The following theorem, in some sense, is known. We do not give the proof
since it exists in the literature, \cite{Conduce}, \cite{Loday}, \cite{Mutlu}%
, \cite{Porter}.

\begin{thm}
The category \textbf{X}$_{2}$\textbf{Mod} of  $2$-crossed modules
is equivalent to the category $\textbf{SimpGrp}_{\leq 2}$ of
simplicial groups with Moore complex of length 2. $\Box $
\end{thm}

\section{\label{s3}Cat$^{2}$-Groups and Crossed Squares}

Although when first introduced by Loday and Walery, \cite{walery},
the notion of crossed square of groups was not linked to that of
cat$^{2}$-groups, it was in this form that Loday gave their
generalisation to an $n$-fold structure, cat$^{n}$-groups (cf.
\cite{Loday}).

\emph{A crossed square} of groups is a commutative square of
groups;
\begin{equation*}
\xymatrix{ L \ar[d]_{\lambda'} \ar[r]^{\lambda} & M \ar[d]^{\mu}
\\ N \ar[r]_{\nu} & P}
\end{equation*}%
together with left actions of $P$ on $L$, $M$, $N$ and a function
$h:M\times N\rightarrow L$. Let $M$ and $N$ act on $M, N$ and $L$
via $P$. The structure must satisfy the following axioms for all
$l\in L$, $m,m^{\prime }\in M$, $n,n^{\prime }\in N$, $p\in P$;
\newline
$(i)$ The homomorphisms $\mu ,\nu ,\lambda ,\lambda ^{\prime }$
and $\mu \lambda $ are crossed modules and both $\lambda ,\lambda ^{\prime }$ are $P$%
-equivariant,\newline
$(ii)$ $h(mm^{\prime },n)=h(^{m}m^{\prime },^{m}n)h(m,n),\newline
(iii)$ $h(m,nn^{\prime })$ $=h(m,n)h(^{n}m,^{n}n^{\prime }),\newline
(iv)$ $\lambda h(m,n)=m^{n}m^{-1},\newline
(v)$ $\lambda ^{\prime }h(m,n)=^{m}nn^{-1},\newline
(vi)$ $h(\lambda l,n)=l^{n}l^{-1},\newline
(vii)$ $h(m,\lambda ^{\prime }l)=^{m}ll^{-1},\newline
(viii)$ $h(^{p}m,^{p}n)=^{p}h(m,n).$

Recall from \cite{Loday} that a cat$^{1}$-group is a triple $(G,s,t)$
consisting of a group $G$ and endomorphisms $s$, the source map, and $t$,
the target map of $G$, satisfying the following axioms:
\begin{equation*}
i)\quad st=t,ts=s,\qquad ii)\quad \lbrack \ker s,\ker t]=1.
\end{equation*}

It was shown that in \cite[Lemma 2.2]{Loday} that setting $C=\ker
s$, $B=\mathrm{im} s$ and $\partial =t|_{C}$, then the conjugation
action makes $\partial :C\rightarrow B$ into a crossed module.
Conversely if $\partial :C\rightarrow B$ is a crossed module, then
setting $G=C\rtimes B$ and letting $s,t$ be defined by
$s(c,b)=(1,b)$ and $t(c,b)=(1,\partial (c)b)$ for $c\in C$, $b\in
B$, then $(G,s,t)$ is a cat$^{1}$-group.

For a cat$^{2}$-group, we again have a group $G$, but this time with two
independent cat$^{1}$-group structures on it. Explicitly:

A cat$^{2}$-group is a 5-tuple, $(G,s_{1},t_{1},s_{2},t_{2})$,
where $(G,s_{i},t_{i})$, $i=1,2$, are cat$^{1}$-groups and
\begin{equation*}
s_{i}s_{j}=s_{j}s_{i},t_{i}t_{j}=t_{j}t_{i},s_{i}t_{j}=t_{j}s_{i}
\end{equation*}%
for $i,j=1,2$, $i\neq j$.

The following proposition was given by Loday (cf. \cite{Loday}).
We only present the sketch proof (see also \cite{mutpor}) of this
result as we need some indication of proofs for later use.

\begin{prop}
\label{loday} \emph{(\cite{Loday})} There is an equivalence of categories
between the category of cat$^{2}$-groups and that of crossed squares.
\end{prop}

\begin{pf}
The cat$^{1}$-group $(G,s_{1},t_{1})$ will give us a crossed
module $\partial :C\rightarrow B$ with $C=\ker s$, $B=\mathrm{im}
s$ and $\partial =t|_{C}$, but as the two cat$^{1}$-group
structures are independent, $(G,s_{2},t_{2})$ restricts to give
cat$^{1}$-group structures on $C$ and $B$ makes $\partial $ a
morphism of cat$^{1}$-groups. We thus get a morphism of crossed
modules
\begin{equation*}
\xymatrix{ \ker s_1\cap\ker s_2 \ar[d] \ar[r] & \mbox{Im}s_1\cap\ker s_2
\ar[d] \\ \ker s_1 \cap \mbox{Im}s_2 \ar[r] & \mbox{Im}s_1\cap\mbox{Im}s_2 }
\end{equation*}%
where each morphism is a crossed module for natural action, i.e.
conjugation in $G$. It remains to produce an $h$-map, but it is
given by the commutator within $G$ since if $x\in \mathrm{im}
s_{1}\cap \ker s_{2}$ and $y\in \ker s_{1}\cap \mathrm{im} s_{2}$
then $[x,y]\in \ker s_{1}\cap \ker s_{2}$. It is easy to check the
crossed square axioms.

Conversely, if
\begin{equation*}
\xymatrix{ L \ar[d] \ar[r] & M \ar[d] \\ N \ar[r] & P }
\end{equation*}
is a crossed square, then we can think of it as a morphism of
crossed modules; $(L,N)\rightarrow (M,P)$.

Using the equivalence between crossed modules and cat$^{1}$-groups this
gives a morphism
\begin{equation*}
\partial :(L\rtimes N,s,t)\longrightarrow (M\rtimes P,s^{\prime },t^{\prime
})
\end{equation*}
of cat$^{1}$-groups. There is an action of $(m,p)\in M\rtimes P$
on $(l,n)\in L\rtimes N$ given by
\begin{equation*}
^{(m,p)}(l,n)=(^{m}(^{p}l)h(m,^{p}n),^{p}n).
\end{equation*}
Using this action, we thus form its associated cat$^{1}$-group
with big group $(L\rtimes N)\rtimes (M\rtimes P)$ and induced
endomorphisms $s_{1},t_{1},s_{2},t_{2}$.
\end{pf}

A generalisation of a crossed square to higher dimensions called a
\textquotedblleft crossed $n$-cube\textquotedblright, was given by
Ellis and Steiner (cf. \cite{Sten}), but we use only the case
$n=2$.

The following result for groups was given by Mutlu and Porter (cf.
\cite{Mutlu}).

Let $\mathbf{G}$ be a simplicial group. Then the following diagram
\begin{equation*}
\xymatrix{ NG_2/\partial_3 NG_3 \ar[d]_{\partial_2'} \ar[r]^-{\partial_2} &
NG_1 \ar[d]^{\mu} \\ \overline{NG}_1 \ar[r]_{\mu'} & G_1 }
\end{equation*}%
is the underlying square of a crossed square. The extra structure
is given as follows;
\newline
$NG_{1}=\ker d_{0}^{1}$ and $\overline{NG}_{1}=\ker d_{1}^{1}$.
Since $G_{1}$ acts on
$NG_{2}/\partial_{3}NG_{3},\overline{NG}_{1}$ and $NG_{1}$, there
are actions of $\overline{NG}_{1}$ on $NG_{2}/\partial _{3}NG_{3}$
and $NG_{1}$ via $\mu ^{\prime }$, and $NG_{1}$ acts on
$NG_{2}/\partial _{3}NG_{3}$ and $\overline{NG}_{1}$ via $\mu $.
Both $\mu $ and $\mu^{\prime }$ are inclusions, and all actions
are given by conjugation. The $h$-map is
\begin{equation*}
\begin{tabular}{crcl}
$h$ : & $NG_{1}\times \overline{NG}_{1}$ & $\longrightarrow $ & $%
NG_{2}/\partial _{3}NG_{3}$ \\
& $(x,\overline{y})$ & $\longmapsto $ & $h(x,y)=[s_{1}x,s_{1}ys_{0}y^{-1}]%
\partial _{3}NG_{3}$.%
\end{tabular}%
\end{equation*}
Here $x$ and $y$ are in $NG_{1}$ as there is a bijection between
$NG_{1}$ and $\overline{NG}_{1}$. This example is clearly
functorial and we denote it by:
$$
\mathbf{M}(-,2) :\textbf{SimpGrp} \longrightarrow
\textbf{Crs$^2$}.
$$
This is the 2-dimensional case of a general construction of a crossed $n$%
-cube from a simplicial group given by Porter, \cite{Porter},
based on some ideas of Loday, \cite{Loday}.

\section{\label{s4}2-Crossed Modules from Crossed Squares}

In this section we will give a description of the passage from
crossed squares to 2-crossed modules by using the `Artin-Mazur'
codiagonal functor and prove directly the $2$-crossed module
structure; a similar construction has been done by Mutlu and
Porter, \cite{mutpor}, in terms of a bisimplicial nerve of a
crossed square.

Conduch{\'{e}} constructed (private communication to Brown in
1984) a $2$-crossed module from a crossed square
\begin{equation*}
\xymatrix{ L \ar[d]_{\lambda'} \ar[r]^{\lambda} & M \ar[d]^{\mu} \\ N
\ar[r]_{\nu} & P }
\end{equation*}%
as
\begin{equation*}
L\llabto{2}{(\lambda ^{-1},\lambda ^{\prime })}M\rtimes N
\llabto{2}{\mu \nu}P.
\end{equation*}
We noted above that the category of crossed modules is equivalent
to that of cat$^{1}$-groups. The corresponding  equivalence in
dimension $2$ is reproved in Proposition \ref{loday}.

We form the associated cat$^{2}$-group. This is
$$
\xymatrix{( L\rtimes N)\rtimes( M\rtimes P)\ar@<-0.5ex>[d]_{s}
\ar@<0.5ex>[d]^{t}\ar@<0.5ex>[r]^-{s'}\ar@<-0.5ex>[r]_-{t'}&M\rtimes
P\ar@<-0.5ex>[d]_{s_{M}}
\ar@<0.5ex>[d]^{t_{M}}\\
N\rtimes P\ar@<0.5ex>[r]^-{s_{N}}\ar@<-0.5ex>[r]_-{t_{N}}&P.}
$$
The source and target maps are defined as follows;
\begin{equation*}
\begin{array}{ll}
s((l,n),(m,p))=(n,p), & s^{\prime }((l,n),(m,p))=(m,p), \\
t((l,n),(m,p))=((\lambda ^{\prime }l)n,\mu (m)p), & t^{\prime
}((l,n),(m,p))=((\lambda l)^{(\nu n)}m,\nu (n)p),\\
s_{N}(n,p)=p,\quad t_{N}(n,p)=\nu (n)p, &s_{M}(m,p)=p,\quad
t_{M}(m,p)=\mu (m)p
\end{array}
\end{equation*}
for $l\in L$, $m\in M$ and $p\in P$.

We take the binerve, that is the nerves in the both directions of
the cat$^{2}$-group constructed. This is a bisimplicial group. The
first few entries in the bisimplicial array are given below
\begin{equation*}
\xymatrix{ \ldots \ar@<-1.5ex>[d]\ar@<-0.5ex>[d] \ar@<0.5ex>[d] \ar@
<1ex>[r] \ar@<0ex>[r] \ar@<-1ex>[r] &
(L^N)\rtimes((L^N)\rtimes(M^P)) \ar@<-1.5ex>[d]\ar@<-0.5ex>[d]
\ar@<0.5ex>[d] \ar@ <0.5ex>[r] \ar@<-0.5ex>[r] & M\rtimes (M\rtimes
P)\ar@<-1.5ex>[d]\ar@<-0.5ex>[d] \ar@<0.5ex>[d]
\\ ((L^L)\rtimes N)\rtimes ((M^M)\rtimes P) \ar@<-0.5ex>[d]
\ar@<0.5ex>[d] \ar@ <1ex>[r] \ar@<0ex>[r] \ar@<-1ex>[r] &
(L^N)\rtimes(M^P) \ar@<-0.5ex>[d] \ar@<0.5ex>[d]\ar@ <0.5ex>[r]
\ar@<-0.5ex>[r] & M\rtimes P\ar@<-0.5ex>[d] \ar@<0.5ex>[d]
\\ N\rtimes (N\rtimes P) \ar@
<1ex>[r] \ar@<0ex>[r] \ar@<-1ex>[r] & (N\rtimes P) \ar@ <1ex>[r]
\ar@<0ex>[r] & P}
\end{equation*}%
where $L^{N}=L\rtimes N$, $M^{P}=M\rtimes P$.

Some reduction has already been done. For example, the double
semi-direct product represents the group of pairs of elements
$((m_{1},p_{1}),(m_{2},p_{2}))\in M\rtimes P$ where $\mu
(m_{1})p_{1}=p_{2}$. This is the group $M\rtimes (M\rtimes P)$,
where the action of $M\rtimes P$ on $M$ is given by
$^{(m,p)}m^{\prime }=$ $^{\mu (m)p}m^{\prime }$.

We will recall the Artin-Mazur codiagonal functor $\nabla $ (cf.
\cite{Artin}) from bisimplicial groups to simplicial groups.

Let \textbf{G.,.} be a bisimplicial group. Put
\begin{equation*}
G_{(n)}=\prod_{p+q=n}G_{p,q}
\end{equation*}
and define $\nabla _{n}\subset G_{(n)}$ as follow; An element $(x_{0},\ldots
,x_{n})$ of $G_{(n)}$ with $x_{p}\in G_{p,n-p}$, is in $\nabla _{n}$ if and
only if
\begin{equation*}
d_{0}^{v}x_{p}=d_{p+1}^{h}x_{p+1}
\end{equation*}%
for each $p=0,\ldots ,n-1$. Next, define the faces and
degeneracies: for $j=0,\ldots ,n$, $D_{j}:\nabla
_{n}\longrightarrow \nabla _{n-1}$ and $\ S_{j}:\nabla
_{n}\longrightarrow \nabla _{n+1}$ by
\begin{equation*}
D_{j}(x)=(d_{j}^{v}x_{0},d_{j-1}^{v}x_{1},\ldots
,d_{1}^{v}x_{j-1},d_{j}^{h}x_{j+1},d_{j}^{h}x_{j+2},\ldots ,d_{j}^{h}x_{n})
\end{equation*}%
\begin{equation*}
S_{j}(x)=(s_{j}^{v}x_{0},s_{j-1}^{v}x_{1},\ldots
,s_{0}^{v}x_{j},s_{j}^{h}x_{j},s_{j}^{h}x_{j+1},\ldots ,s_{j}^{h}x_{n}).
\end{equation*}%
Thus $\nabla (\mathbf{G.,.})=\{\nabla _{n}:D_{j},S_{j}\}$ is a simplicial
group.

We now examine this construction in low dimension:\newline
EXAMPLE:

For $n=0$, $G_{(0)}=G_{0,0}$. For $n=1$, we have
\begin{equation*}
\nabla _{1}\subset G_{(1)}=G_{1,0}\times G_{0,1}
\end{equation*}%
where
\begin{equation*}
\nabla _{1}=\{(g_{1,0},g_{0,1}):d_{0}^{v}(g_{1,0})=d_{1}^{h}(g_{0,1})\}
\end{equation*}%
together with the homomorphisms
\begin{align*}
D_{0}^{1}(g_{1,0},g_{0,1})&=(d_{0}^{v}g_{1,0},d_{0}^{h}g_{0,1}),\\
D_{1}^{1}(g_{1,0},g_{0,1})&=(d_{1}^{v}g_{1,0},d_{1}^{h}g_{0,1}),\\
S_{0}^{0}(g_{0,0})&=(s_{0}^{v}g_{0,0},s_{0}^{h}g_{0,0}).
\end{align*}
For $n=2$, we have
\begin{equation*}
\nabla _{2}\subset G_{(2)}=\prod_{p+q=2}G_{p,q}=G_{2,0}\times G_{1,1}\times
G_{0,2}
\end{equation*}
where
\begin{equation*}
\nabla
_{2}=\{(g_{2,0},g_{1,1},g_{0,2}):d_{0}^{v}(g_{2,0})=d_{1}^{h}(g_{1,1}),\text{
}d_{0}^{v}(g_{1,1})=d_{2}^{h}(g_{0,2})\}.
\end{equation*}
Now, we use the Artin-Mazur codiagonal functor to obtain a simplicial group $%
\mathbf{G}$ (of some complexity).

The base group is still $G_{0}\cong P$. However the group of
1-simplices  is the subset of
\begin{equation*}
G_{1,0}\times G_{0,1}=(M\rtimes P)\times (N\rtimes P),
\end{equation*}
consisting of $(g_{1,0},g_{0,1})=((m,p),(n,p^{\prime }))$ where $\mu
(m)p=p^{\prime }$, i.e.,
\begin{equation*}
G_{1}=\{((m,p),(n,p^{\prime })):d_{0}^{v}(m,p)=\mu (m)p=p^{\prime
}=d_{1}^{h}(n,p^{\prime })\}.
\end{equation*}
We see that the composite of two elements
\begin{center}
$(m_{1},p_{1},n_{1},\mu (m_{1})p_{1})$ and $(m_{2},p_{2},n_{2},\mu
(m_{2})p_{2})$
\end{center}
 becomes
\begin{equation*}
(m_{1}\ ^{p_{1}}m_{2},p_{1}p_{2},n_{1}\text{ }^{\mu (m_{1})p_{1}}n_{2},\mu
(m_{1}\ ^{p_{1}}m_{2})p_{1}p_{2})
\end{equation*}%
(by the inter-change law). The subgroup $G_{1}$ of these elements is
isomorphic to $N\rtimes (M\rtimes P),$ where $M$ acts on $N$ via
$P$, $ ^{m}n=$ $ ^{\mu m}n$. Indeed, one can easily show that the
map
\begin{equation*}
\begin{array}{ccll}
f: & G_{1} & \longrightarrow & N\rtimes (M\rtimes P) \\
& (m,p,n,\mu (m)p) & \longmapsto & (n,m,p)%
\end{array}%
\end{equation*}%
is an isomorphism.

Identifying $G_{1}$ with $N\rtimes (M\rtimes P)$, $d_{0}$ and
$d_{1}$ have the descriptions
\begin{align*}
d_{0}(n,m,p) & =  \upsilon (n)\mu (m)p \\
d_{1}(n,m,p) & =  p.
\end{align*}
We next turn to the group of 2-simplices: this is the subset
$G_{2}$ of
\begin{equation*}
G_{2,0}\times G_{1,1}\times G_{0,2}=M\rtimes (M\rtimes P)\times ((L\rtimes
N)\rtimes (M\rtimes P))\times (N\rtimes (N\rtimes P))
\end{equation*}
whose elements
\begin{equation*}
((m_{2},m_{1},p),(l,n,m,p^{\prime }),(n_{2},n_{1},p^{\prime \prime }))
\end{equation*}
are such that
\begin{align*}
d_{0}^{v}(m_{2},m_{1},p) & =  d_{1}^{h}(l,n,m,p^{\prime }) \\
d_{0}^{v}(l,n,m,p^{\prime }) & = d_{2}^{h}(n_{2},n_{1},p^{\prime \prime }).%
\end{align*}
This gives the relations between the individual coordinates implying
that $m_{1}=m$, $\mu (m_{2})p=p^{\prime }$, $\lambda ^{\prime }(l)n=n_{2}$ and $%
\mu (m)p^{\prime }=p^{\prime \prime }$. Thus the elements of $G_{2}$
have the form
\begin{equation*}
((m_{2},m_{1},p),((l,n),(m_{1},\mu \left( m_{2}\right) p)),(\lambda ^{\prime
}(l)n,n_{1},\mu (m_{1}m_{2})p)).
\end{equation*}
We then deduce the isomorphism
\begin{equation*}
f:G_{2}\longrightarrow (L\rtimes (N\rtimes M))\rtimes (N\rtimes
(M\rtimes P))
\end{equation*}%
given by \begin{multline*} ((m_{2},m_{1},p),((l,n),(m_{1},\mu
\left( m_{2}\right) p)),((\lambda ^{\prime }l)n,n_{1},\mu
(m_{1}m_{2})p))\\ \longmapsto \left(
(l,(n,m_{1})),(n_{1},(m_{2},p))\right) .
\end{multline*}
Therefore we can get a 2-truncated simplicial group
$\mathbf{G}^{(2)}$ that looks like
\begin{equation*}
\xymatrix@C=3pc{\textbf{G}^{(2)}: (L \rtimes (N \rtimes M))\rtimes
(N \rtimes (M \rtimes P))\ar@<2ex>[r]^-{d_{0}^{2} ,d_{1}^{2},
d_{2}^{2}}\ar@<1ex>[r]\ar[r] &  N \rtimes (M\rtimes
P)\ar@<1ex>[r]^-{d_{0}^{1}
,d_{1}^{1}}\ar[r]\ar@<1ex>[l]\ar@<2ex>[l]^-{s_{0}^{1}, s_{1}^{1}} &
P \ar@<1ex>[l]^-{s_{0}^{0}} }
\end{equation*}
with the faces and degeneracies;
\begin{equation*}
d_{0}^{1}(n,m,p)=\nu (n)\mu (m)p,\quad d_{1}^{1}(n,m,p)=p,\quad
s_{0}(p)=(1,1,p)
\end{equation*}
and \begin{align*} d_{0}^{2}((l,(n,m_{1})),(n_{1},(m_{2},p))) & =
(n_{1},(\lambda l)^{\nu
(n)}m_{1},\nu (n)\mu (m_{2})p), \\
d_{1}^{2}((l,(n,m_{1})),(n_{1},(m_{2},p))) & =  (n_{1}(\lambda
^{\prime
}l)n,m_{1}m_{2},p), \\
d_{2}^{2}((l,(n,m_{1})),(n_{1},(m_{2},p))) & =  (n,m_{2},p), \\
s_{0}^{1}(n,m,p) & =  ((1,(1,m)),(n,(1,p))), \\
s_{1}^{1}(n,m,p) & =  ((1,(n,1)),(1,(m,p))).%
\end{align*}

For the verification of the simplicial identities, see appendix.

\textbf{Remark: }

The construction given above may be shortened in terms of the $\overline{W}$
construction or `bar' construction (cf. \cite{Artin}, \cite{bullejos}), but
we have not attempted this method.

Loday, \cite{Loday}, defined the mapping cone of a complex as
analogous to the construction of the Moore complex of a simplicial
group. (for further work see also \cite{Con1}). We next describe
the mapping cone of a crossed square of groups as follows:

\begin{prop}
\label{pr1} The Moore complex of the simplicial group
$\mathbf{G}^{(2)}$ is the mapping cone, in the sense of Loday, of
the crossed square. Furthermore, this mapping cone complex has a
2-crossed module structure of groups.
\end{prop}

\begin{pf}
Given the $2$-truncated simplicial group $\mathbf{G}^{(2)}$
described above, look at its Moore complex; we have
$NG_{0}=G_{0}=P$. The second term of the Moore complex is
$NG_{1}=\ker d_{0}^{1}$. By the definition of $d_{0}^{1}$,
$(n,m,p)\in \ker d_{0}^{1}$ if and only if $p=\mu (m)^{-1}\nu
(n)^{-1}$. Since $d_{0}^{1}(n^{-1},m^{-1},\mu (m)\nu (n))=\nu
(n)^{-1}\mu (m)^{-1}\mu (m)\nu (n)=1,$ we have $(n^{-1},m^{-1},\mu
(m)\nu (n))\in \ker d_{0}^{1}$. Furthermore there is an isomorphism
$f_{1}: NG_{1}\longrightarrow M\rtimes N$ given by
\begin{equation*}
(n^{-1},m^{-1},\mu (m)\nu (n))\mapsto (m,n).
\end{equation*}
We note that via this isomorphism, the map $\partial _{1}:M\rtimes
N\rightarrow P$ is given by $\partial _{1}(m,n)=\mu (m)\nu (n)$.

Now we investigate the intersection of the kernels of $d_{0}^{2}$
and $d_{1}^{2}$. Let
\begin{equation*}
\mathbf{x}=((l,(n,m_{1})),(n_{1},(m_{2},p)))\in (L\rtimes (N\rtimes
M))\rtimes (N\rtimes (M\rtimes P)).
\end{equation*}
If $\mathbf{x}\in \ker d_{0}^{2}$, by the definition of
$d_{0}^{2}$, we have
\begin{equation*}
n_{1}=1,\text{ }(\lambda l)^{\nu (n)}m_{1}=1,\text{ }\nu (n)\mu (m_{2})p=1.
\end{equation*}%
If $\mathbf{x}\in \ker d_{1}^{2}$, by the definition of
$d_{1}^{2}$, we have
\begin{equation*}
n_{1}(\lambda ^{\prime }l)n=1,\text{ }m_{1}m_{2}=1,\text{ }p=1.
\end{equation*}
>From these equalities we have $n=(\lambda ^{\prime }l)^{-1}$, and
from
\begin{align*}1&=(\lambda l)^{\nu (n)}m_{1}\\
&=(\lambda l)^{\nu ((\lambda ^{\prime }l)^{-1})}m_{1}\\
&=(\lambda l)^{\mu \lambda l^{-1}}m_{1}  \qquad (\mu \lambda =\nu
\lambda ^{\prime })\\
&=(\lambda l)(\lambda l)^{-1}m_{1}(\lambda l)\\
&= m_{1}(\lambda l),
\end{align*}
we have $m_{1}=m_{2}^{-1}=(\lambda l)^{-1}$ and $p=1$. Therefore, $\mathbf{x%
}\in \ker d_{0}^{2}\cap \ker d_{1}^{2}$ if and only if
\begin{equation*}
\mathbf{x=(}(l,(\lambda ^{\prime }l^{-1},\lambda l^{-1})),(1,\lambda l,1)).
\end{equation*}%
Thus we get $\ker d_{0}^{2}\cap \ker d_{1}^{2}$ $\cong L$.

>From these calculations, we have
\begin{align*}
d_{2}|_{\ker d_{0}^{2}\cap \ker d_{1}^{2}}((l,(\lambda ^{\prime
}l^{-1},\lambda l^{-1})),(1,\lambda l,1))&=(\lambda ^{\prime
}l^{-1},\lambda l,1).
\end{align*}%
Of course $(\lambda ^{\prime }l^{-1},\lambda l,1)\in NG_{1}$ since
\begin{equation*}
d_{0}^{1}(\lambda ^{\prime }l^{-1},\lambda l,1)=\nu \lambda ^{\prime
}l^{-1}\mu (\lambda l)1=1.
\end{equation*}%
By using above isomorphism $f_{1}$ and $d_{2}|_{\ker d_{0}^{2}\cap
\ker d_{1}^{2}}$, we can identify the map $\partial _{2}$ on $L$
by
\begin{align*}
\partial _{2}(l) & =  f_{1}d_{2}|_{\ker d_{0}^{2}\cap \ker
d_{1}^{2}}((l,(\lambda ^{\prime }l^{-1},\lambda l^{-1})),(1,\lambda l,1)) \\
& =  f_{1}(\lambda ^{\prime }l^{-1},\lambda l,1) \\
& =  (\lambda l^{-1},\lambda ^{\prime }l)\in M\rtimes N.%
\end{align*}
It can be seen that $\partial _{2}$ and $\partial _{1}$ are homomorphisms
and
\begin{align*}
\partial _{1}\partial _{2}(l) & =  \partial _{1}(\lambda l^{-1},\lambda
^{\prime }l) \\
& =  \mu (\lambda l)\nu (\lambda ^{\prime }l^{-1}) \\
& =  1\quad (\text{by }\nu \lambda ^{\prime }=\mu \lambda ).%
\end{align*}%
Thus, if given a crossed square
\begin{equation*}
\xymatrix{ L \ar[d]_{\lambda'} \ar[r]^{\lambda} & M \ar[d]^{\mu} \\ N
\ar[r]_{\nu} & P }
\end{equation*}%
its mapping cone complex is
$$
L\llabto{2}{\partial_2}M\rtimes N\llabto{2}{\partial_1}P
$$
where $\partial _{2}l=(\lambda l^{-1},\lambda ^{\prime }l)$ and
$\partial _{1}(m,n)=\mu (m)\nu (n).$ The semi-direct product
$M\rtimes N$ can be formed by making $N$ acts on $M$ via $P$,
$^{n}m=$ $^{\nu (n)}m$, where the $P$-action is the given one.

These elementary calculations are useful as they pave the way for the
calculation of the Peiffer commutator of $x=(m,n)$ and $y=(c,a)$ in the
above complex;
\begin{align*}
\left\langle x,y\right\rangle &={}^{\partial _{1}x}yxy^{-1}x^{-1} \\
&={}^{\mu (m)\nu (n)}(c,a)(m,n)(^{a^{-1}}{c^{-1}},a^{-1})(^{n^{-1}}{m^{-1}}%
,n^{-1}) \\
&={}(^{\mu (m)\nu (n)}c,^{\mu (m)\nu (n)}a)(m^{\nu
(na^{-1})}(c^{-1}){^{\nu (n^{-1}n^{-1})}}m^{-1},na^{-1}n^{-1})
\end{align*}%
which on multiplying out and simplifying is
$$
(^{\nu (nan^{-1})}mm^{-1},{}^{\mu (m)}(nan^{-1})(na^{-1}n^{-1}))
$$
(Note that any dependence on $c$ vanishes!)

Conduch{\'{e} (unpublished work)} defined the Peiffer lifting for
this structure by
\begin{equation*}
\{x,y\}=\{(m,n),(c,a)\}=h(m,nan^{-1}).
\end{equation*}%
For the axioms of 2-crossed module see appendix.
\end{pf}

We thus have two ways of going from simplicial groups to 2-crossed
modules
\newline
$(i)$ (\cite{Mutlu}) directly to get
\begin{equation*}
NG_{2}/\partial _{3}NG_{3}\longrightarrow NG_{1}\longrightarrow NG_{0},
\end{equation*}%
$(ii)$ indirectly via the square axiom $\mathbf{M}(\mathbf{G}, 2)$
and then by the above construction to get
\begin{equation*}
NG_{2}/\partial _{3}NG_{3}\longrightarrow \ker d_{0}\rtimes \ker
d_{1}\longrightarrow G_{1},
\end{equation*}%
and they clearly give the same homotopy type. More precisely
$G_{1}$ decomposes as $\ker d_{1}\rtimes s_{0}G_{0}$ and the $\ker
d_{0}$ factor in the middle term of $(ii)$ maps down to that in
this decomposition by the identity map. Thus $d_{0}$ induces a
quotient map from $(ii)$ to $(i)$ with
kernel isomorphic to%
\begin{equation*}
1\longrightarrow \ker d_{0}\overset{=}{\longrightarrow }\ker d_{0}
\end{equation*}%
which is thus contractible.

\textbf{Note}: The construction given above from a crossed square to
a 2-crossed module preserves the homotopy type. In fact, Ellis (cf.
\cite{ellis}) defined the homotopy groups of the crossed square is
the homology groups of
the complex%
$$
L\llabto{2}{\partial_2}M\rtimes
N\llabto{2}{\partial_1}P\longrightarrow 1
$$
where $\partial _{1}$ and $\partial _{2}$ are defined above.

\section{\label{s5}Quadratic Modules from 2-Crossed Modules}

Quadratic modules of groups were initially defined by Baues,
\cite{baus1, baus2}, as models for connected 3-types. In this
section we will define a functor from the category
$\textbf{X$_2$Mod}$\textbf{X}$_{2}$\textbf{Mod} of 2- crossed
modules to that of quadratic modules $\textbf{QM}$. Before giving
the definition of quadratic module we should recall some
structures.

Recall that \textit{a pre-crossed module} is a group homomorphism
$\partial :M\rightarrow N$ together with an action of $N$ on $M$,
written $^{n}m$ for $n\in N$ and $m\in M$, satisfying the
condition $\partial (^{n}m)=n\partial (m)n^{-1}$ for all $m\in M$
and $n\in N$.

\textit{A nil(2)-module} is a pre-crossed module $\partial
:M\rightarrow N$ with an additional \textquotedblleft
nilpotency\textquotedblright condition. This condition is
$P_{3}(\partial )=1$, where $P_{3}(\partial )$ is the subgroup of
$M$ generated by Peiffer commutator $\left\langle
x_{1},x_{2},x_{3}\right\rangle $ of length $3$.

The Peiffer commutator in a pre-crossed module $\partial :M\rightarrow N$ is
defined by
\begin{equation*}
\left\langle x,y\right\rangle =(^{\partial x}y)xy^{-1}x^{-1}
\end{equation*}%
for $x,y\in M$.

For a group $G$, the group
\begin{equation*}
G^{ab}=G/[G,G]
\end{equation*}%
is the abelianization of $G$ and
\begin{equation*}
\partial ^{cr}:M^{cr}=M/P_{2}(\partial )\rightarrow N
\end{equation*}%
is the crossed module associated to the pre-crossed module
$\partial :M\rightarrow N$. Here $P_{2}(\partial )=\langle
M,M\rangle$ is the Peiffer subgroup of $M$.\newline The following
definition is due to Baues (cf. \cite{baus2}).

\begin{defn}
\emph{A quadratic module} $(\omega ,\delta ,\partial )$ is a diagram
\begin{equation*}
\xymatrix{ & C\otimes C\ar[dl]_{\omega} \ar[d]^{w} \\ L \ar[r]_{\delta} & M
\ar[r]_{\partial} & N }
\end{equation*}%
of homomorphisms between groups such that the following axioms are satisfied.%
\newline
$\mathbf{QM1)}$ The homomorphism $\partial :M\rightarrow N$ is a
nil(2)-module with Peiffer commutator map $w$ defined above. The quotient
map $M\twoheadrightarrow C=(M^{cr})^{ab}$ is given by $x\mapsto \overline{x}%
, $ where $\overline{x}\in C$ denotes the class represented by $x\in M$ and $%
C=(M^{cr})^{ab}$ is the abelianization of the associated crossed module $%
M^{cr}\rightarrow N$.\newline
$\mathbf{QM2)}$ The boundary homomorphisms $\partial $ and $\delta $ satisfy
$\partial \delta =1$ and the quadratic map $\omega $ is a lift of the
Peiffer commutator map $w$, that is $\delta \omega =w$ or equivalently
\begin{equation*}
\delta \omega (\overline{x}\otimes \overline{y})=(^{\partial
x}y)xy^{-1}x^{-1}=\left\langle x,y\right\rangle
\end{equation*}%
for $x,y\in M$.\newline
$\mathbf{QM3)}$\textbf{\ }$L$ is an $N$-group and all homomorphisms of the
diagram are equivariant with respect to the action of $N$. Moreover, the
action of $N$ on $L$ satisfies the formula ($a\in L,x\in M$)
\begin{equation*}
^{\partial x}a=\omega (\left( \overline{x}\otimes \overline{\delta a}\right)
\left( \overline{\delta a}\otimes \overline{x}\right) )a.
\end{equation*}%
$\mathbf{QM4)}$ Commutators in $L$ satisfy the formula ($a,b\in L$)
\begin{equation*}
\omega (\overline{\delta a}\otimes \overline{\delta b})=[b,a].
\end{equation*}
\end{defn}

A map $\varphi :(\omega ,\delta ,\partial )\rightarrow (\omega ^{\prime
},\delta ^{\prime },\partial ^{\prime })$ between quadratic modules is given
by a commutative diagram, $\varphi =(l,m,n)$
\begin{equation*}
\xymatrix{ C\otimes C\ar[d]_{\varphi_{\ast}\otimes\varphi_{\ast}}
\ar[r]^-{\omega} & L \ar[d]_{l} \ar[r]^{\delta} & M \ar[d]_{m}
\ar[r]^{\partial} & N \ar[d]_{n} \\ C^{\prime }\otimes C^{\prime }
\ar[r]_-{\omega^{\prime}} & L^{\prime} \ar[r]_{\delta^{\prime}} & M^{\prime}
\ar[r]_{\partial^{\prime}} & N^{\prime} }
\end{equation*}%
where $(m,n)$ is a morphism between pre-crossed modules which induces $%
\varphi _{\ast }:C\rightarrow C^{\prime }$ and where $l$ is an $n$%
-equivariant homomorphism. Let $\textbf{QM}$ be the category of
quadratic modules and of maps as in above diagram.

Now, we construct a functor from the category of 2-crossed modules to that
of quadratic modules.

Let
$$
C_2\llabto{2}{\partial_2}C_1 \llabto{2}{\partial_1}C_0
$$
be a 2-crossed module. Let $P_{3}$ be the subgroup of $C_{1}$
generated by elements of the form
\begin{equation*}
\left\langle \left\langle x,y\right\rangle ,z\right\rangle \text{ and }%
\left\langle x,\left\langle y,z\right\rangle \right\rangle
\end{equation*}
with $x,y,z\in C_{1}$. We obtain $\partial _{1}(\left\langle
\left\langle x,y\right\rangle ,z\right\rangle )=1$ and $\partial
_{1}(\left\langle x,\left\langle y,z\right\rangle \right\rangle
)=1$, since $\partial _{1}$ is a pre-crossed module.

Let $P_{3}^{\prime }$ be the subgroup of $C_{2}$ generated by elements of
the form
\begin{equation*}
\{\left\langle x,y\right\rangle ,z\}\text{ and }\{x,\left\langle
y,z\right\rangle \}
\end{equation*}
for $x,y,z\in C_{1},$ where $\left\{ -,-\right\} $ is the Peiffer lifting
map. Then there are quotient groups
\begin{equation*}
M=C_{1}/P_{3}
\end{equation*}
and
\begin{equation*}
L=C_{2}/P_{3}^{\prime }.
\end{equation*}
Then, $\partial : M\rightarrow C_{0}$ given by $\partial
(xP_{3})=\partial_{1}(x)$ is a well defined group homomorphism
since $\partial _{1}(P_{3})=1$. We thus get the following
commutative diagram
\begin{equation*}
\xymatrix{ C_1\ar[dr]_{q_1}\ar[rr]^{\partial_1}& & C_0\\
&M\ar[ur]_{\partial} }
\end{equation*}
where $q_{1}:C_{1}\rightarrow M$ is the quotient map.

Furthermore, from the first axiom of $2$-crossed module
$\mathbf{2CM1)}$, we can write $\partial _{2}\{\left\langle
x,y\right\rangle ,z\}=\left\langle \left\langle x,y\right\rangle
,z\right\rangle $ and $\partial _{2}\{x,\left\langle
y,z\right\rangle \}=\left\langle x,\left\langle y,z\right\rangle
\right\rangle .$ Therefore, the map $\delta :L\rightarrow M$ given
by $\delta (lP_{3}^{\prime })=(\partial _{2}l)P_{3}$ is a well
defined group homomorphism since $\partial _{2}(P_{3}^{\prime
})=P_{3}$.

Thus we get the following commutative diagram;
\begin{equation*}
\xymatrix{ & C\otimes C\ar[dl]_{\omega} \ar[d]^{w} \\ L \ar[r]^{\delta } & M
\ar[r]^{\partial } & N \ar@{=}[d] \\ C_2 \ar[u]_{q_2} \ar[r]_{\partial_2} &
C_1 \ar[u]_{q_1} \ar[r]_{\partial_1} & C_0 }
\end{equation*}%
where $q_{1}$ and $q_{2}$ are the quotient maps and $C=(M^{cr})^{ab}$ is a
quotient of $C_{1}$. The quadratic map is given by the Peiffer lifting map
\begin{equation*}
\left\{ -,-\right\} :C_{1}\times C_{1}\longrightarrow C_{2},
\end{equation*}%
namely
\begin{equation*}
\omega (\overline{x^{\prime }}\otimes \overline{y^{\prime
}})=q_{2}(\{x,y\})
\end{equation*}%
for $x^{\prime },y^{\prime }\in M$ and $x,y\in C_{1}$.

\begin{prop}
\label{pr2} The diagram%
\begin{equation*}
\xymatrix{ & C\otimes C\ar[dl]_{\omega} \ar[d]^{w} \\ L \ar[r]_{\delta } & M
\ar[r]_{\partial } & N }
\end{equation*}%
is a quadratic module of groups.
\end{prop}

\begin{pf}
For the axioms, see appendix.
\end{pf}

\begin{prop}
\label{ho2}The homotopy groups of the 2-crossed module are isomorphic to
that of its associated quadratic module.
\end{prop}

\begin{pf}
Consider the 2-crossed module
\begin{equation}
\label{1}C_2\llabto{2}{\partial_2}C_1 \llabto{2}{\partial_1}C_0
\tag{1}
\end{equation}
and its associated quadratic module
\begin{equation}
\label{2}\xymatrix{ & C\otimes C\ar[dl]_{\omega} \ar[d]^{w} \\ L
\ar[r]_{\delta } & M \ar[r]_-{\partial } &N= C_0.}  \tag{2}
\end{equation}%
The homotopy groups of (\ref{1}) are
\begin{equation*}
\pi_{i}=
\begin{cases}
C_{0}/\partial _{1}(C_{1}) &   i=1, \\
\ker \partial _{1}/\mathrm{im} \partial _{2} &   i=2, \\
\ker \partial _{2} &   i=3, \\
0 &   i=0 \text{ or }i>3.%
\end{cases}%
\end{equation*}%
The homotopy groups of (\ref{2}) are
\begin{equation*}
\pi _{i}^{\prime }= \begin{cases}
C_{0}/\partial (M) &   i=1, \\
\ker \partial /\mathrm{im} \delta  &   i=2, \\
\ker \delta  &   i=3, \\
0 &   i=0 \text{ or } i>3.%
\end{cases}
\end{equation*}%

We claim that $\pi _{i}=\pi _{i}^{\prime }$ for all $i\geq 0.$ In
fact, since $\partial (M)\cong \partial _{1}(C_{1})$, clearly $\pi
_{1}\cong \pi_{1}^{\prime }$. Also $\ker \partial =\dfrac{\ker
\partial _{1}}{P_{3}}$, $\mathrm{im} \delta \cong \dfrac{\mathrm{im} \partial
_{2}}{P_{3}}$ so that $\pi_{2}^{\prime }=\dfrac{\ker \partial
_{1}/P_{3}}{\mathrm{im} \partial _{2}/P_{3}} \cong \dfrac{\ker
\partial _{1}}{\mathrm{im} \partial _{2}}\cong $ $\pi _{2}$. Consider now
$\pi _{3}^{\prime }=\{xP_{3}^{\prime }:\partial _{2}(x)\in
P_{3}\}$. We show that given $xP_{3}^{\prime }\in \pi _{3}^{\prime
}$, there is $x^{\prime }P_{3}^{\prime }\in \pi _{3}^{\prime }$
with $xP_{3}^{\prime }=x^{\prime }P_{3}^{\prime }$ and $x^{\prime
}\in \ker \partial _{2}$. In fact, observe that since $\partial
_{2}\{\left\langle x,y\right\rangle ,z\}=\left\langle \left\langle
x,y\right\rangle ,z\right\rangle ,$ $\partial _{2}\{x,\left\langle
y,z\right\rangle \}=\left\langle x,\left\langle y,z\right\rangle
\right\rangle$, we have $\partial _{2}(P_{3}^{\prime })=P_{3}.$
Hence $\partial _{2}(x)\in P_{3}$ implies $\partial
_{2}(x)=\partial _{2}(w),$ $w\in P_{3}^{\prime };$ thus $\partial
_{2}(xw^{-1})=1;$ then take $x^{\prime }=xw^{-1},$ so that
$xP_{3}^{\prime
}=x^{\prime }P_{3}^{\prime }$ and $\partial _{2}(x^{\prime })=1.$ Define $%
\alpha :\pi _{3}^{\prime }\rightarrow \pi _{3},$ $\alpha (xP_{3}^{\prime
})=\alpha (x^{\prime }P_{3}^{\prime })=x^{\prime }$ and $\beta :\pi
_{3}\rightarrow \pi _{3}^{\prime },$ $\beta (x)=xP_{3}.$ Clearly $\alpha $
and $\beta $ are inverse bijections, proving the claim. It follows that $%
(1)$ and $(2)$ represent the same homotopy type.
\end{pf}

\section{\label{s6}Simplicial Groups and Quadratic Modules}

Baues gives a construction of a quadratic module from a simplicial
group in Appendix B to Chapter IV of \cite{baus2}. The quadratic
modules can be given by using higher dimensional Peiffer elements in
verifying the axioms.

This section is a brief r{\'{e}}sum{\'{e}}  defining a variant of
the Carrasco-Cegarra pairing operators that are called \textit{Peiffer Pairings }%
(cf. \cite{cc}). The construction depends on a variety of sources,
mainly Conduch\'{e}, \cite{Conduce},  Mutlu and Porter,
\cite{Mutlu, mutlu2,mutpor}. We define a normal subgroup $N_{n}$
of $G_{n}$ and a set $P(n)$ consisting of pairs of elements
$(\alpha ,\beta )$ from $S(n)$ (cf. \cite{Mutlu}) with $\alpha
\cap \beta =\emptyset $ and $\beta <\alpha $ , with respect to the
lexicographic ordering in $S(n)$ where $\alpha =(i_{r},\ldots,
i_{1}),\beta =(j_{s},\ldots, j_{1})\in S(n)$. The pairings that we
will need,
\begin{equation*}
\{F_{\alpha ,\beta }:NG_{n-\sharp \alpha }\times NG_{n-\sharp \beta
}\rightarrow NG_{n}:(\alpha ,\beta )\in P(n),n\geq 0\}
\end{equation*}%
are given as composites by the diagram%
\begin{equation*}
\xymatrix{ NG_{n-\sharp\alpha } \times NG_{n-\sharp\beta
}\ar[d]_{s_{\alpha}\times s_{\beta}} \ar[r]^-{F_{\alpha,\beta}} & NG_n \\
G_n \times G_n \ar[r]_{\mu} & G_n \ar[u]_{p} }
\end{equation*}%
where $s_{\alpha }=s_{i_{r}},\ldots, s_{i_{1}}:NG_{n-\sharp \alpha
}\rightarrow G_{n},$\quad $s_{\beta }=s_{j_{s}},\ldots,
s_{j_{1}}:NG_{n-\sharp \beta }\rightarrow G_{n},$
$p:G_{n}\rightarrow NG_{n}$ is defined by composite projections
$p(x)=p_{n-1}\ldots p_{0}(x)$, where $p_{j}(z)=zs_{j}d_{j}(z)^{-1}$
with $j=0,1,\ldots,n-1$ and $\mu :G_{n}\times G_{n}\rightarrow
G_{n}$ is given by the commutator map and $\sharp \alpha $  is the
number of the elements in the set of $\alpha$;  similarly for
$\sharp \beta $. Thus
\begin{equation*}
F_{\alpha ,\beta }(x_{\alpha },y_{\beta })=p[s_{\alpha }(x_{\alpha
}),s_{\beta }(x_{\beta })].
\end{equation*}

\begin{defn}
Let $N_{n}$ or more exactly $N_{n}^{G}$ be the normal subgroup of $G_{n}$
generated by elements of the form $F_{\alpha ,\beta }(x_{\alpha },y_{\beta
}) $ where $x_{\alpha }\in NG_{n-\sharp \alpha }$ and $y_{\beta }\in
NG_{n-\sharp \beta }$ .
\end{defn}

This normal subgroup $N_{n}^{G}$ depends functorially on $G$, but
we will usually abbreviate $N_{n}^{G}$ to $N_{n},$ when no change
of group is involved. Mutlu and Porter (cf. \cite{Mutlu})
illustrate this normal subgroup for $n=2,3,4$, but we only
consider for $n=3$.
\begin{ex}
For all $x_{1}\in NG_{1},y_{2}\in NG_{2},$ the corresponding generators of $%
N_{3}$ are: \begin{align*} F_{(1,0)(2)}(x_{1},y_{2})
&=[s_{1}s_{0}x_{1},s_{2}y_{2}][s_{2}y_{2},s_{2}s_{0}x_{1}], \\
F_{(2,0)(1)}(x_{1},y_{2})
&=[s_{2}s_{0}x_{1},s_{1}y_{2}][s_{1}y_{2},s_{2}s_{1}x_{1}][s_{2}s_{1}x_{1},s_{2}y_{2}][s_{2}y_{2},s_{2}s_{0}x_{1}]\\
\intertext{ and for all $x_{2}\in NG_{2},y_{1}\in NG_{1},$}
F_{(0)(2,1)}(x_{2},y_{1})&=[s_{0}x_{2},s_{2}s_{1}y_{1}][s_{2}s_{1}y_{1},s_{1}x_{2}][s_{2}x_{2},s_{2}s_{1}y_{1}]\\
\intertext{
whilst for all $x_{2},y_{2}\in NG_{2},$}%
F_{(0)(1)}(x_{2},y_{2})
&=[s_{0}x_{2},s_{1}y_{2}][s_{1}y_{2},s_{1}x_{2}][s_{2}x_{2},s_{2}y_{2}], \\
F_{(0)(2)}(x_{2},y_{2}) &=[s_{0}x_{2},s_{2}y_{2}], \\
F_{(1)(2)}(x_{2},y_{2})
&=[s_{1}x_{2},s_{2}y_{2}][s_{2}y_{2},s_{2}x_{2}].
\end{align*}
\end{ex}

The following theorem is proved by Mutlu and Porter (cf.
\cite{mutlu2}).

\begin{thm}
Let $\mathbf{G}$ be a simplicial group and for $n>1$, let $D_{n}$
the subgroup of $G_{n}$ generated by degenerate elements. Let
$N_{n}$ be the normal subgroup generated by elements of the form
$F_{\alpha ,\beta }(x_{\alpha },y_{\beta })$ with $(\alpha ,\beta
)\in P(n)$ where $x_{\alpha }\in NG_{n-\sharp \alpha }$ and
$y_{\beta }\in NG_{n-\sharp \beta }$. Then
\begin{equation*}
NG_{n}\cap D_{n}=N_{n}\cap D_{n}
\end{equation*}%
\begin{flushright}
$\Box $
\end{flushright}
\end{thm}
Baues defined a functor from the category of simplicial groups to
that of quadratic modules (cf. \cite{baus2}). Now we will
reconstruct this functor by using  the $F_{\alpha ,\beta }$
functions. We will use the $F_{\alpha ,\beta }$ functions in
verifying the axioms of quadratic module.

Let $\mathbf{G}$ be a simplicial group with Moore complex
$\mathbf{NG}$. Suppose that $G_{3}=D_{3}$. Notice  that
$P_{3}(\partial _{1})$ is the subgroup of $NG_{1}$ generated by
triple brackets
\begin{equation*}
\left\langle x,\left\langle y,z\right\rangle \right\rangle \text{ and }%
\left\langle \left\langle x,y\right\rangle ,z\right\rangle
\end{equation*}%
for $x,y,z\in NG_{1}$. Let $P_{3}^{\prime }(\partial _{1})$ be the
subgroup of $NG_{2}/\partial _{3}NG_{3}$ generated by elements of
the form
\begin{align*}
\omega (\left\langle x,y\right\rangle ,z)&=s_{0}(\left\langle
x,y\right\rangle )s_{1}zs_{0}(\left\langle x,y\right\rangle
)^{-1}s_{1}(\left\langle x,y\right\rangle
)s_{1}z^{-1}s_{1}(\left\langle x,y\right\rangle )^{-1}\\
\intertext{and} \omega (x,\left\langle y,z\right\rangle
)&=s_{0}xs_{1}(\left\langle y,z\right\rangle
)s_{0}x^{-1}s_{1}xs_{1}(\left\langle y,z\right\rangle
)^{-1}s_{1}x^{-1}.
\end{align*}
Then we have quotient groups
\begin{align*}
M&=NG_{1}/P_{3}(\partial _{1})\\
\intertext{and} L&=(NG_{2}/\partial _{3}NG_{3})/P_{3}^{\prime
}(\partial _{1}).
\end{align*}
We obtain $\overline{\partial _{2}}\omega (\left\langle x,y\right\rangle
,z)=\left\langle \left\langle x,y\right\rangle ,z\right\rangle $ and $%
\overline{\partial _{2}}\omega (x,\left\langle y,z\right\rangle
)=\left\langle x,\left\langle y,z\right\rangle \right\rangle .$ Thus, $%
\delta :L\rightarrow M$ given by $\delta (\overline{a}P_{3}^{\prime
}(\partial _{1}))=\overline{\partial _{2}}(\overline{a})P_{3}(\partial _{1})$
is a well defined group homomorphism, where $\overline{a}$ is a coset in $%
NG_{2}/\partial _{3}NG_{3}$.

Therefore, we obtain the following diagram,
\begin{equation*}
\xymatrix{ & C\otimes C\ar[dl]_{\omega} \ar[d]^{w} \\ L \ar[r]^{\delta } & M
\ar[r]^{\partial } & N \ar@{=}[d] \\ NG_2/\partial_3 NG_3 \ar[u]_{q_2}
\ar[r]_-{\overline{\partial_2}} & NG_1 \ar[u]_{q_1} \ar[r]_{\partial_1} &
NG_0 }
\end{equation*}%
where $q_{1}$ and $q_{2}$ are quotient maps and $\delta q_{2}=q_{1}\overline{%
\partial _{2}},\partial q_{1}=\partial _{1},$ and the quadratic map $\omega $
is defined by
\begin{equation*}
\omega (\{q_{1}x\}\otimes \{q_{1}y\})=q_{2}(\overline{%
s_{0}xs_{1}ys_{0}x^{-1}s_{1}xs_{1}y^{-1}s_{1}x^{-1}})
\end{equation*}%
for $x,y\in NG_{1},$ $q_{1}x,q_{1}y\in M$ and $\{q_{1}x\}\otimes
\{q_{1}y\}\in C\otimes C$ and where $C=((M)^{cr})^{ab}$.

\begin{prop}
\label{pr3} The diagram%
\begin{equation*}
\xymatrix{ & C\otimes C\ar[dl]_{\omega} \ar[d]^{w} \\ L \ar[r]_{\delta } & M
\ar[r]_{\partial } & N }
\end{equation*}%
is a quadratic module of groups.
\end{prop}

\begin{pf}
We show that all the axioms of quadratic module are verified by
using the functions $F_{\alpha ,\beta }$ in the appendix.
\end{pf}

Alternatively, this proposition can be reproved differently, by
making use of the 2-crossed module constructed from a simplicial
group by Mutlu and Porter (cf. \cite{Mutlu}). We now give a sketch
of the argument. In \cite{Mutlu}, it is shown that given a
simplicial group $\mathbf{G}$, one can construct a 2-crossed
module

\begin{equation}
\label{3}\xymatrix{ NG_2/\partial_3(NG_3 \cap D_3)
\ar[r]^-{\overline{\partial_2}} & NG_1 \ar[r]^{\partial_1 } & NG_0
}  \tag{3}
\end{equation}%
where $\{x,y\}=s_{0}xs_{1}ys_{0}x^{-1}s_{1}y^{-1}s_{1}x^{-1}\partial
_{3}(NG_{3}\cap D_{3})$ for $x,y\in NG_{1}.$ \newline
Clearly we have a commutative diagram
\begin{equation*}
\xymatrix{ NG_2/\partial_3(NG_3 \cap D_3)\ar[d]^{j}
\ar[r]^-{\overline{\partial_2}} & NG_1 \ar[r]^{\partial_1 }\ar@{=}[d] & NG_0
\ar@{=}[d]\\ NG_2/\partial_3(NG_3)\ar[r]&NG_1 \ar[r]^{\partial_1}&NG_0.}
\end{equation*}%
Consider now the quadratic module associated to the 2-crossed
module (\ref{3}), as in Section \ref{s5} of this paper.
\begin{equation*}
\xymatrix{ & C\otimes C\ar[dl]_{\omega'} \ar[d] \\ L' \ar[r]^{\delta' } & M
\ar[r] & N \ar@{=}[d] \\ \dfrac{{NG_2}}{{\partial_3 (NG_3\cap D_3)}}
\ar[u]_{q'_2} \ar[r]_-{\overline{\partial_2}} & NG_1 \ar[u]
\ar[r]_{\partial_1} & NG_0 }
\end{equation*}%
Then one can see that $L^{\prime}=\Omega /\partial _{3}(NG_{3}\cap D_{3}),$ where $%
\Omega $ is the subgroup of $NG_{2}$ generated by elements of the
form $$ s_{0}(\left\langle x,y\right\rangle
)s_{1}zs_{0}(\left\langle x,y\right\rangle
)^{-1}s_{1}(\left\langle x,y\right\rangle
)s_{1}z^{-1}s_{1}(\left\langle x,y\right\rangle )^{-1} $$ and

$$ s_{0}xs_{1}(\left\langle y,z\right\rangle
)s_{0}x^{-1}s_{1}xs_{1}(\left\langle y,z\right\rangle
)^{-1}s_{1}x^{-1}.
$$%
On the other hand we have, from Section \ref{s5}, $L=\Omega
/\partial _{3}(NG_{3})$. Hence there is a map $i:L^{\prime
}\rightarrow L$ with
\begin{equation}
\label{4}\omega =i\omega ^{\prime },\quad \delta ^{\prime }=\delta
i. \tag{4}
\end{equation}%
Since
\begin{equation*}
\xymatrix{ C\otimes C \ar[r]^-{\omega'}&L'\ar[r]^{\delta'} &M
\ar[r]^{\partial } & N }
\end{equation*}%
is, by construction a quadratic module, it is straightforward to
check, using (\ref{4}),  that
\begin{equation*}
\xymatrix{ C\otimes C \ar[r]^-{\omega}&L\ar[r]^{\delta} &M
\ar[r]^{\partial } & N }
\end{equation*}%
is also a quadratic module.

\begin{prop}
Let $\mathbf{G}$ be a simplicial group, let $\pi'_{i}$ be the
homotopy groups of its associated quadratic module and let $\pi_i$
be the homotopy groups of the classifying space of $\mathbf{G}$;
then $\pi_i \cong\pi'_i$ for $i=0,1,2,3$.
\end{prop}

\begin{pf}
Let $\mathbf{G}$ be a simplicial group. The $n$th homotopy groups of
$\mathbf{G}$ is the $n$th homology of the Moore complex of
$\mathbf{G}$, i.e.,
\begin{equation*}
\pi _{n}(\mathbf{G})\cong H_{n}(\mathbf{NG})\cong \frac{\ker
d_{n-1}^{n-1}\cap NG_{n-1}}{d_{n}^{n}(NG_{n})}.
\end{equation*}%
Thus the homotopy groups $\pi _{n}(\mathbf{G})=\pi _{n}$ of $\mathbf{G}$ are%
\begin{equation*}
\pi _{n}= \begin{cases}
NG_{0}/d_{1}(NG_{1}) &  n=1, \\
\\
\dfrac{{\ker d_{1}\cap NG_{1}}}{{d_{2}(NG_{2})}}  & n=2, \\
\\
\dfrac{\ker d_{2}\cap NG_2}{d_3({NG_3})}   & n=3, \\
\\
0 &   n=0 \text{ or } n>3.%
\end{cases}
\end{equation*}%
and the homotopy groups $\pi _{n}^{\prime }$ of its associated quadratic
module are
\begin{equation*}
\pi _{n}^{\prime }=
\begin{cases}
NG_{0}/\partial (M)   & n=1, \\
\ker \partial /\mathrm{im} \delta   & n=2, \\
\ker \delta &   n=3, \\
0 &   n=0 \text{ or } n>3.%
\end{cases}%
\end{equation*}%
We claim that $\pi _{n}^{\prime }\cong \pi _{n}$ for $n=1,2,3$.
Since $M=NG_{1}/P_{3}(\partial _{1})$ and $d_{1}(P_{3}(\partial
_{1}))=1$, we have
\begin{equation*}
\partial (M)=\partial (NG_{1}/P_{3}(\partial _{1}))=d_{1}(NG_{1})
\end{equation*}%
and then
\begin{equation*}
\pi _{1}^{\prime }=NG_{0}/\partial (M)\cong NG_{0}/d_{1}(NG_{1})=
\pi _{1}.
\end{equation*}%
Also $\ker \partial = \dfrac{\ker d_{1}\cap NG_1}{P_{3}(\partial
_{1})}$ and $\mathrm{im} \delta=d_{2}(NG_{2})/P_{3}(\partial
_{1})$ so that we have
\begin{equation*}
\pi _{2}^{\prime }=\frac{\ker \partial }{\mathrm{im} \delta
}=\frac{(\ker
d_{1}\cap NG_1)/P_{3}(\partial _{1})}{d_{2}(NG_{2})/P_{3}(\partial _{1})}\cong \frac{%
\ker d_{1}\cap NG_1}{d_{2}(NG_{2})}= \pi _{2}.
\end{equation*}%
We know that $P_{3}^{\prime }(\partial _{1})$ is generated by elements of
the form
$$
s_{0}(\left\langle x,y\right\rangle )s_{1}zs_{0}(\left\langle
x,y\right\rangle )^{-1}s_{1}(\left\langle x,y\right\rangle
)s_{1}z^{-1}s_{1}(\left\langle x,y\right\rangle )^{-1}
$$
and
$$%
s_{0}xs_{1}(\left\langle y,z\right\rangle
)s_{0}x^{-1}s_{1}xs_{1}(\left\langle y,z\right\rangle )^{-1}s_{1}x^{-1}.
$$
Since \begin{multline*} d_{2}(s_{0}(\left\langle x,y\right\rangle
)s_{1}zs_{0}(\left\langle x,y\right\rangle
)^{-1}s_{1}(\left\langle x,y\right\rangle
)s_{1}z^{-1}s_{1}(\left\langle x,y\right\rangle )^{-1})\\
\begin{aligned}
& = {}^{d_{1}\left\langle x,y\right\rangle }z\left\langle
x,y\right\rangle
(z^{-1})\left\langle x,y\right\rangle ^{-1} \\
 &=  \left\langle \left\langle x,y\right\rangle ,z\right\rangle
\in P_{3}(\partial _{1})\end{aligned}
\end{multline*}%
and \begin{multline*} d_{2}(s_{0}xs_{1}(\left\langle
y,z\right\rangle )s_{0}x^{-1}s_{1}xs_{1}(\left\langle
y,z\right\rangle )^{-1}s_{1}x^{-1})\\ \begin{aligned}& =
s_{0}d_{1}x\left\langle y,z\right\rangle
s_{0}d_{1}x^{-1}(x)\left\langle
y,z\right\rangle ^{-1}x^{-1} \\
& = {} ^{d_{1}x}\left\langle y,z\right\rangle x\left\langle
y,z\right\rangle
^{-1}x^{-1} \\
& =  \left\langle x,\left\langle y,z\right\rangle \right\rangle
\in P_{3}(\partial _{1})\end{aligned}%
\end{multline*}
we get $d_{2}(P_{3}^{\prime }(\partial _{1}))=P_{3}(\partial _{1})$.
The isomorphism between $\pi _{3}^{\prime }$ and $\pi _{3}$ can be
proved similarly  to the proof of Proposition \ref{ho2}.
\end{pf}

\section{\label{s7}Quadratic Modules from Crossed Squares}

In this section, we will define a functor from the category of crossed
squares to that of quadratic modules. Our construction can be briefly
explained as:

Given a crossed square, we consider the associated 2-crossed
module ( from Section \ref{s4}) and then we build the quadratic
module corresponding to this 2-crossed module (from Section
\ref{s5}). In other words, we are just composing two functors.
Thus, there is no need to worry in this section about direct
proofs, as they hold automatically from the results of Sections
\ref{s4} and \ref{s5}. In particular, the homotopy type is clearly
preserved, as it is preserved at each step.

Now let
\begin{equation*}
\xymatrix{ L\ar[r]^{\lambda}\ar[d]_{\lambda'}& M\ar[d]^{\mu}\\ N
\ar[r]_{\nu}&P}
\end{equation*}%
be a crossed square of groups. Consider its associated 2-crossed
module from Section \ref{s4}%
\begin{equation*}
L\llabto{2}{(\lambda ^{-1},\lambda ^{\prime })}M\rtimes N%
\llabto{2}{\mu \nu }P.
\end{equation*}%
>From this 2-crossed module, we can get a quadratic module as in
Section \ref{s5}%
\begin{equation*}
\xymatrix{ & C\otimes C\ar[dl]_{\omega} \ar[d]^{w} \\ C_2 \ar[r]_-{\delta} &
C_1 \ar[r]_-{\partial} & C_0 }
\end{equation*}%
where $C_{0}=P,$ $C_{1}=(M\rtimes N)/P_{3},$ $C_{2}=L/P_{3}^{\prime }$, $%
C=((C_{1})^{cr})^{ab}$ and the quadratic map is given by
\begin{equation*}
\begin{array}{cccc}
\omega : & C\otimes C & \longrightarrow & C_{2} \\
& [q_{1}(m,n)]\otimes \lbrack q_{1}(c,a)] & \longmapsto &
q_{2}(h(m,nan^{-1}))%
\end{array}%
\end{equation*}%
for $(m,n),(c,a)$ $\in M\rtimes N$ , $q_{1}(m,n),$ $q_{1}(c,a)\in C_{1}$ and
$[q_{1}(m,n)]\otimes \lbrack q_{1}(c,a)]\in C\otimes C$. Furthermore $P_{3}$
is the subgroup of $M\rtimes N$ generated by elements of the form
\begin{equation*}
\left\langle \left\langle (m,n),(c,a)\right\rangle ,(m^{\prime },n^{\prime
})\right\rangle \text{ and }\left\langle (m,n),\left\langle (c,a),(m^{\prime
},n^{\prime })\right\rangle \right\rangle
\end{equation*}%
for $(m,n),(c,a),(m^{\prime },n^{\prime })\in M\rtimes N,$ and $%
P_{3}^{\prime }$ is the subgroup of $L$ generated by elements of the form
\begin{equation*}
h\left( ^{\nu (nan^{-1})}mm^{-1},^{\upsilon (^{\mu
(m)}(nan^{-1})(na^{-1}n^{-1}))}n^{\prime }\right)
\end{equation*}%
and
\begin{equation*}
h\left( m,^{\nu (n)}\left( ^{\mu (c)}(an^{\prime }a^{-1})(an^{\prime
}{}^{-1}a^{-1})\right) \right)
\end{equation*}%
for $(m,n),(c,a),(m^{\prime },n^{\prime })$ $\in M\rtimes N.$
$\delta :C_{2}\longrightarrow C_{1}$ is defined by $\delta
(lP_{3}^{\prime })=(\lambda l^{-1},\lambda ^{\prime }l)P_{3}$ and
$\partial :C_{1}\rightarrow C_{0}$ is defined by $\partial
(q_{1}(m,n))=\mu (m)\upsilon (n)$.

The proof of the axioms of quadratic module is similar to the proof of the
axioms of Proposition \ref{pr1}.

\section{Appendix}

 \textbf{The proof of simplicial identities}:%
\begin{align*}
\begin{split}
d_{0}^{2}s_{0}^{1}(n,m,p)
&=d_{0}^{2}((1,(1,m)),(n,(1,p)))  \\
& =(n,\lambda 1.^{\nu (1)}m,\upsilon (1)\mu (1)p)\\
& =(n,m,p)=id,\end{split}
\begin{split}
d_{1}^{2}s_{0}^{1}(n,m,p)
&=d_{1}^{2}((1,(1,m)),(n,(1,p)))\\
&=(n\lambda ^{\prime }1.1,m,p)\\
&=(n,m,p)=id
\end{split}
\end{align*}
and%
\begin{align*}
\begin{split}
d_{1}^{2}s_{1}^{1}(n,m,p)
&= d_{1}^{2}((1,(n,1)),(1,(m,p)))  \\
& =(1.\lambda ^{\prime }1.n,1.m,p)\\
& =(n,m,p)=id,\end{split}
\begin{split}
d_{2}^{2}s_{1}^{1}(n,m,p)
&=d_{2}^{2}((1,(n,1)),(1,(m,p)))\\
&=(n,m,p)=id\\
&
\end{split}
\end{align*}
and%
\begin{equation*}
d_{2}^{2}s_{0}^{1}(n,m,p)=d_{2}((1,(1,m)),(n,(1,p)))=(1,1,p)=s_{0}^{0}(p)=s_{0}^{0}d_{1}^{1}(n,m,p).
\end{equation*}%
Similarly%
\begin{equation*}
d_{1}^{1}d_{1}^{2}((l,(n,m_{1})),(n_{1},(m_{2},p)))=d_{1}^{1}(n_{1}(\lambda
^{\prime }l)n,m_{1}m_{2},p)=p
\end{equation*}%
and
\begin{equation*}
d_{1}^{1}d_{2}^{2}((l,(n,m_{1})),(n_{1},(m_{2},p)))=d_{1}^{1}(n,m_{2},p)=p
\end{equation*}%
then we have $d_{1}^{1}d_{1}^{2}=d_{1}^{1}d_{2}^{2}.$%
\begin{align*}
d_{0}^{1}d_{1}^{2}((l,(n,m_{1})),(n_{1},(m_{2},p))) & =
d_{0}^{1}(n_{1}(\lambda ^{\prime }l)n,m_{1}m_{2},p) \\
& =  \nu (n_{1})\nu \lambda ^{\prime }(l)\nu (n)\mu (m_{1})\mu (m_{2})p \\
\\
d_{0}^{1}d_{0}^{2}((l,(n,m_{1})),(n_{1},(m_{2},p))) & =
d_{0}^{1}(n_{1},(\lambda l)^{\nu (n)}m_{1},\nu (n)\mu (m_{2})p) \\
& =  \nu (n_{1})\mu ((\lambda l)^{\nu (n)}m_{1})\nu (n)\mu (m_{2})p \\
& =  \nu (n_{1})\mu (\lambda l)\nu (n)\mu (m_{1})\nu (n)^{-1}\nu
(n)\mu (m_{2})p \\
& =  \nu (n_{1})\nu \lambda ^{\prime }(l)\nu (n)\mu (m_{1})\mu (m_{2})p%
\text{ \qquad }(\nu \lambda ^{\prime }=\mu \lambda )%
\end{align*}%
so $d_{0}^{1}d_{1}^{2}=d_{0}^{1}d_{0}^{2}.$
\begin{flushright}
$\Box $
\end{flushright}

\textbf{The Proof of Axioms (Proposition \ref{pr1}):}

$\mathbf{2CM1)}$
\begin{align*}
\partial _{2}\{x,y\} &=(\lambda h(m,nan^{-1})^{-1},\lambda ^{\prime
}h(m,nan^{-1})) \\
&=(^{\nu (nan^{-1})}mm^{-1},^{\mu (m)}(nan^{-1})(na^{-1}n^{-1})) \\
&=\left\langle x,y\right\rangle
\end{align*}%
by axioms of the crossed square.

$\mathbf{2CM2)}$ We will show that $\{\partial
_{2}(l_{0}),\partial _{2}(l_{1})\}=[l_{1},l_{0}].$ As $\partial
_{2}l=(\lambda l^{-1},\lambda ^{\prime }l)$, this need the
calculation of $h(\lambda l_{0}^{-1},\lambda
^{\prime }(l_{0}l_{1}l_{0}^{-1})); $ but the crossed square axioms $%
h(\lambda l,n)=l^{n}l^{-1}$ and $h(m,\lambda ^{\prime }l)=(^{m}l)l^{-1}$
together with the fact that the map $\lambda :L\rightarrow M$ is a crossed
module, give:
\begin{align*}
h(\lambda l_{0}^{-1},\lambda ^{\prime }(l_{0}l_{1}l_{0}^{-1}))
&={}^{\mu
\lambda (l_{0})^{-1}}(l_{0}l_{1}l_{0}^{-1}).l_{0}l_{1}^{-1}l_{0}^{-1} \\
&= [l_{1},l_{0}].
\end{align*}

$\mathbf{2CM3)}$ For the elements of $M\rtimes N$ are; $m=(m_{0},n_{0}),m^{%
\prime }=(m_{1},n_{1}),m^{\prime \prime }=(m_{2},n_{2})$ we have

 $(i)$ %
\begin{align*}
\{mm^{\prime },m^{\prime \prime }\}
&=\{(m_{0},n_{0})(m_{1},n_{1}),(m_{2},n_{2})\} \\
&=\{(m_{0}^{\nu (n_{0})}(m_{1}),n_{0}n_{1}),(m_{2},n_{2})\} \\
&=h(m_{0}^{\nu (n_{0})}(m_{1}),n_{0}n_{1}n_{2}n_{1}^{-1}n_{0}^{-1}) \\
&=h(^{\mu (m_{0})\nu (n_{0})}m_{1},^{\mu (m_{0})\nu
(n_{0})}(n_{1}n_{2}n_{1}^{-1}))h(m_{0},n_{0}n_{1}n_{2}n_{1}^{-1}n_{0}^{-1})
\\
&={}^{\mu (m_{0})\nu
(n_{0})}h(m_{1},n_{1}n_{2}n_{1}^{-1})h(m_{0},n_{0}n_{1}n_{2}n_{1}^{-1}n_{0}^{-1})
\\
&={}^{\partial _{1}(m)}(\{m^{\prime },m^{\prime \prime
}\})h(m_{0},n_{0}n_{1}n_{2}n_{1}^{-1}n_{0}^{-1}).
\end{align*}%
Since
\begin{align*}
mm^{\prime \prime }m^{\prime ^{-1}}
&=(m_{1},n_{1})(m_{2},n_{2})(m_{1},n_{1})^{-1} \\
&=(m_{1},n_{1})(m_{2},n_{2})(^{\nu (n_{1}^{-1})}(m_{1}^{-1}),n_{1}^{-1}) \\
&=(m_{1}^{\nu (n_{1})}m_{2},n_{1}n_{2})(^{\nu
(n_{1}^{-1})}(m_{1}^{-1}),n_{1}^{-1}) \\
&=(m_{1}^{\nu (n_{1})}m_{2}{^{\nu (n_{1}n_{2}n_{1}^{-1})}}%
(m_{1})^{-1},n_{1}n_{2}n_{1}^{-1})
\end{align*}%
and
\begin{align*}
\{m,m^{\prime }m^{\prime \prime }m^{\prime ^{-1}}\}
&=\{(m_{0},n_{0}),(m_{1}^{\nu (n_{1})}m_{2}{^{\nu (n_{1}n_{2}n_{1}^{-1})}}%
(m_{1})^{-1},n_{1}n_{2}n_{1}^{-1})\} \\
&=h(m_{0},n_{0}n_{1}n_{2}n_{1}^{-1}n_{0}^{-1}),
\end{align*}%
we get%
\begin{equation*}
^{\partial _{1}(m)}(\{m^{\prime },m^{\prime \prime
}\})h(m_{0},n_{0}n_{1}n_{2}n_{1}^{-1}n_{0}^{-1})={}^{\partial
_{1}(m)}(\{m^{\prime },m^{\prime \prime }\})\{m,m^{\prime
}m^{\prime \prime }m^{\prime ^{-1}}\},
\end{equation*}%
and thus
\begin{equation*}
\{mm^{\prime },m^{\prime \prime }\}={}^{\partial
_{1}(m)}(\{m^{\prime },m^{\prime \prime }\})\{m,m^{\prime
}m^{\prime \prime }m^{\prime ^{-1}}\}.
\end{equation*}
$(ii)$ \begin{align*} \{m,m^{\prime }m^{\prime \prime }\}
&=\{(m_{0},n_{0}),(m_{1}^{\nu
(n_{1})}m_{2},n_{1}n_{2})\} \\
&=h(m_{0},n_{0}n_{1}n_{2}n_{0}^{-1}) \\
&=h(m_{0},n_{0}n_{1}n_{0}^{-1}n_{0}n_{2}n_{0}^{-1}) \\
&=h(m_{0},n_{0}n_{1}n_{0}^{-1})h(^{\nu
(n_{0}n_{1}n_{0}^{-1})}m_{0},^{n_{0}n_{1}n_{0}^{-1}}n_{0}n_{2}n_{0}^{-1}) \\
&=\{m,m^{\prime }\}h(^{\nu
(n_{0}n_{1}n_{0}^{-1})}m_{0},^{n_{0}n_{1}n_{0}^{-1}}n_{0}n_{2}n_{0}^{-1})
\end{align*}%
and this gives the following result
\begin{equation*}
\{m,m^{\prime }m^{\prime \prime }\}=\{m,m^{\prime }\}^{mm^{\prime
}(m^{-1})}\{m,m^{\prime \prime }\}.
\end{equation*}

$\mathbf{2CM4)}$
\begin{align*}
\{x,\partial _{2}l\}\{\partial _{2}l,x\} & =  \{(m,n),(\lambda
l^{-1},\lambda ^{\prime }l)\}\{(\lambda l^{-1},\lambda ^{\prime
}l),(m,n)\}
\\
& =  h(m,n\lambda ^{\prime }ln^{-1})h(\lambda l^{-1},\lambda
^{\prime
}ln\lambda ^{\prime }l^{-1}) \\
& =  h(m,\lambda ^{\prime }(^{n}l))h(\lambda (l^{-1}),\lambda
^{\prime
}ln\lambda ^{\prime }l^{-1}) \\
& = {} ^{\mu (m)\nu (n)}l^{\nu (n)}(l^{-1})(l^{-1}){^{\nu \lambda
^{\prime
}l\nu (n)\nu \lambda ^{\prime }l^{-1}}}l%
\end{align*}
and this simplifies as expected to give the correct result.
\begin{flushright}
$\Box $
\end{flushright}

\textbf{The Proof of Axioms (Proposition \ref{pr2}):}

$\mathbf{QM1)}$ Clearly $\partial :M\rightarrow N$ is a
nil(2)-module as the Peiffer commutators which in the forms
$\left\langle x,\left\langle y,z\right\rangle \right\rangle $ and
$\left\langle \left\langle x,y\right\rangle ,z\right\rangle $ are
in $P_{3}(\partial _{1}).$

$\mathbf{QM2)}$ It is easy to see that $\delta \partial =1.$ Also
\begin{align*}
\delta \omega (\overline{x^{\prime }}\otimes \overline{y^{\prime }})
&=\delta q_{2}(\{x,y\}) \\
&=q_{1}\partial _{2}\{x,y\} \\
&=q_{1}(^{\partial _{1}x}y)x(y)^{-1}(x)^{-1} \\
&=(^{\partial x^{\prime }}y^{\prime })x^{\prime }(y^{\prime
})^{-1}(x^{\prime })^{-1}.
\end{align*}%
for $\overline{x^{\prime }},\overline{y^{\prime }}\in C$ and $x^{\prime
},y^{\prime }\in M.$

$\mathbf{QM3)}$ For $x^{\prime }\in M$ and $[a]\in L,$
\begin{align*}
\omega \left( \overline{x^{\prime }}\otimes \overline{(\partial
_{2}[a])^{\prime }}\overline{(\partial _{2}[a])^{\prime }}\otimes \overline{%
x^{\prime }}\right) [a] &=q_{2}((\{x,\partial _{2}a\}\{\partial
_{2}a,x\})a)\qquad \text{ (by definition)} \\
&=q_{2}(^{\partial _{1}x}a(^{x}a^{-1})(^{x}a)a^{-1}a) \qquad(\text{by}\quad \mathbf{2CM4)}) \\
&={}^{\partial x^{\prime }}[a].
\end{align*}

$\mathbf{QM4)}$
\begin{align*}
\omega \left( \overline{\delta \lbrack a]}\otimes \overline{\delta \lbrack b]%
}\right) &=\omega \left( \overline{(\partial _{2}a)^{\prime
}}\otimes
\overline{(\partial _{2}b)^{\prime }}\right) \qquad \text{ (by commutativity)%
} \\
&=q_{2}\{\partial _{2}a,\partial _{2}b\} \\
&=[[b],[a]].
\end{align*}

for $[a],[b]\in L$.
\begin{flushright}
$\Box $
\end{flushright}

\textbf{The Proof of Axioms (Proposition \ref{pr3}):}

We display the elements omitting the overlines in our calculation
to save complication.

$\mathbf{QM1)}$  $\partial :M\rightarrow N$ is a nil(2)-module as
the Peiffer commutators which in the forms $\left\langle
x,\left\langle y,z\right\rangle \right\rangle $ and $\left\langle
\left\langle x,y\right\rangle ,z\right\rangle $ are in
$P_{3}(\partial _{1}).$

$\mathbf{QM2)}$ For all $q_{1}x,q_{1}y\in M,$%
\begin{align*}
\delta \omega (\{q_{1}x\}\otimes \{q_{1}y\}) &=\delta
q_{2}(s_{0}xs_{1}ys_{0}x^{-1}s_{1}xs_{1}y^{-1}s_{1}x^{-1}) \\
&=q_{1}\partial _{2}(s_{0}xs_{1}ys_{0}x^{-1}s_{1}xs_{1}y^{-1}s_{1}x^{-1}) \\
&=q_{1}(d_{2}(s_{0}xs_{1}ys_{0}x^{-1}s_{1}xs_{1}y^{-1}s_{1}x^{-1})) \\
&=q_{1}(s_{0}d_{1}xys_{0}d_{1}x^{-1}xy^{-1}x^{-1}) \\
&=\left\langle q_{1}x,q_{1}y\right\rangle \qquad \text{ by
}\partial q_{1}=\partial _{1}.
\end{align*}

$\mathbf{QM3)}$ Supposing $D_{3}=G_{3},$ we know from \cite{Mutlu}
that
\begin{equation*}
d_{3}(F_{(2,0)(1)}(x,a))=[s_{0}x,s_{1}d_{2}a][s_{1}d_{2}a,s_{1}x][s_{1}x,a][a,s_{0}x]\in
\partial _{3}(NG_{3}).
\end{equation*}%
>From this equality we have%
\begin{equation*}
\lbrack s_{0}x,s_{1}d_{2}a][s_{1}d_{2}a,s_{1}x]\equiv \lbrack
s_{0}x,a][a,s_{1}x] \mod \partial _{3}(NG_{3}).
\end{equation*}%
Thus we get
\begin{align*}
\omega (\{q_{1}x\}\otimes \{\delta q_{2}a\}) &=\omega \left(
\{q_{1}x\}\otimes \{q_{1}\partial _{2}a\}\right) \\
&=q_{2}\left(
s_{0}(x)s_{1}d_{2}(a)s_{0}(x)^{-1}s_{1}(x)s_{1}d_{2}(a)^{-1}s_{1}(x)^{-1}%
\right) \\
&\equiv q_{2}([s_{0}x,a][a,s_{1}x]) \\
&={}^{\partial q_{1}(x)}q_{2}(a)q_{2}(^{x}(a^{-1})),
\end{align*}%
and similarly from
\begin{equation*}
d_{3}(F_{(0)(2,1)}(a,x))=[s_{0}d_{2}a,s_{1}x][s_{1}x,s_{1}d_{2}a][a,s_{1}x]%
\in \partial _{3}(NG_{3}\cap D_{3})=\partial _{3}(NG_{3}),
\end{equation*}%
we have
\begin{align*}
\omega (\{\delta q_{2}a\}\otimes \{q_{1}x\}) &=\omega
(\{q_{1}\partial
_{2}a\}\otimes \{q_{1}x\}) \\
&=q_{2}(s_{0}d_{2}(a)s_{1}(x)s_{0}d_{2}(a)^{-1}s_{1}d_{2}(a)s_{1}(x)^{-1}s_{1}d_{2}(a)^{-1})
\\
&\equiv q_{2}([s_{1}x,a]) \\
&=q_{2}(^{x}a)q_{2}a^{-1}.
\end{align*}%
Consequently we have,
\begin{equation*}
\omega \left( \{q_{1}x\}\otimes \{\delta q_{2}a\}\{\delta
q_{2}a\}\otimes \{q_{1}x\}\right) q_{2}a={}^{\partial
q_{1}(x)}q_{2}a.
\end{equation*}%

$\mathbf{QM4)}$  From \cite{Mutlu}, we get
\begin{equation*}
d_{3}(F_{(0),(1)}(a,b))=[s_{0}d_{2}a,s_{1}d_{2}b][s_{1}d_{2}b,s_{1}d_{2}a][a,b]\in \partial _{3}(NG_{3}\cap D_{3})=\partial _{3}(NG_{3}).
\end{equation*}
>From this equality, we can write
\begin{equation*}
\lbrack s_{0}d_{2}a,s_{1}d_{2}b][s_{1}d_{2}b,s_{1}d_{2}a]\equiv
\lbrack b,a] \mod \partial _{3}(NG_{3}).
\end{equation*}%
Thus we have
\begin{align*}
\omega (\{\delta q_{2}a\}\otimes \{\delta q_{2}b\}) &=\omega
(\{q_{1}\partial _{2}a\}\otimes \{q_{1}\partial _{2}b\}) \\
&=q_{2}(s_{0}d_{2}(a)s_{1}d_{2}(b)s_{0}d_{2}(a)^{-1}s_{1}d_{2}(a)s_{1}d_{2}(b)^{-1}s_{1}d_{2}(a)^{-1})
\\
&\equiv [q_{2}b,q_{2}a].
\end{align*}%
for $q_{2}a,q_{2}b\in L.$
\begin{flushright}
$\Box $
\end{flushright}

\bigskip

\end{document}